\documentclass[10pt]{amsart}
\usepackage{amssymb}
\usepackage{a4}
\usepackage[latin1]{inputenc}
\usepackage[T1]{fontenc}
\theoremstyle{plain}
\newtheorem{theorem}[subsection]{Theorem}
\newtheorem{corollary}[subsection]{Corollary}
\newtheorem{lemma}[subsection]{Lemma}
\newtheorem{proposition}[subsection]{Proposition}

\theoremstyle{definition}

\theoremstyle{remark}
\newtheorem{remark}[subsection]{Remark}
\numberwithin{equation}{section}
\newcommand{\ko}{{\mathcal O}}
\newcommand{\kh}{{\mathcal H}}
\newcommand{\IC}{{\mathbb C}}
\newcommand{\IF}{{\mathbb F}}
\newcommand{\IN}{{\mathbb N}}
\newcommand{\IP}{{\mathbb P}}
\newcommand{\IQ}{{\mathbb Q}}

\newcommand{\IZ}{{\mathbb Z}}
\newcommand{\gotha}{{\mathfrak a}}
\newcommand{\gothg}{{\mathfrak g}}
\newcommand{\gothh}{{\mathfrak h}}
\newcommand{\gothm}{{\mathfrak m}}
\newcommand{\gothp}{{\mathfrak p}}

\newcommand{\lra}{\longrightarrow}
\newcommand{\xra}{\xrightarrow}
\newcommand{\xla}{\xleftarrow}
\newcommand{\isom}{\cong}
\newcommand{\tensor}{\otimes}

\DeclareMathOperator{\Hom}{Hom}
\DeclareMathOperator{\End}{End}
\DeclareMathOperator{\Ext}{Ext}
\DeclareMathOperator{\Tor}{Tor}
\DeclareMathOperator{\PAut}{PAut}
\DeclareMathOperator{\Aut}{Aut}
\DeclareMathOperator{\Kern}{Ker}

\DeclareMathOperator{\Liegl}{\mathfrak{gl}}
\DeclareMathOperator{\LieGl}{Gl}
\DeclareMathOperator{\LiePGl}{PGl}

\DeclareMathOperator{\Quot}{Quot}
\DeclareMathOperator{\Hilb}{Hilb}
\DeclareMathOperator{\Pic}{Pic}

\DeclareMathOperator{\sing}{sing}
\DeclareMathOperator{\reg}{reg}

\DeclareMathOperator{\ch}{ch}
\DeclareMathOperator{\tr}{tr}

\DeclareMathOperator{\even}{even}
\DeclareMathOperator{\Spec}{Spec}

\DeclareMathOperator{\rk}{rk}
\DeclareMathOperator{\td}{td}
\DeclareMathOperator{\gr}{gr}
\DeclareMathOperator{\ini}{in}
\DeclareMathOperator{\NS}{NS}
\DeclareMathOperator{\codim}{codim}

\newcommand{\GIT}{/\!\!/}

\newcommand{\dual}{{\scriptstyle\vee}}
\newcommand{\tra}{^t}

\begin{document}
\title{Singular symplectic moduli spaces}
\author{D.~Kaledin, M.~Lehn, and Ch.~Sorger}
%\date{1.~April 2005, 19:00 Uhr}

\begin{abstract} Moduli spaces of semistable sheaves on a K3 or
abelian surface with respect to a general ample divisor are shown to be
locally factorial, with the exception of symmetric products of a
K3 or abelian surface and the class of moduli spaces found
by O'Grady. Consequently, since singular moduli space that do not belong to these
exceptional cases have singularities in codimension $\geq4$ they
do no admit projective symplectic resolutions.
\end{abstract}

\address{
Dmitry Kaledin\\
Independent University of Moscow\\
B. Vlassievski per. 11\\
Moscow, 119002, Russia}\email{kaledin@mccme.ru}

\thanks{The first author has been partially supported by
CRDF Award RM1-2354-MO02 and the
``cooperation franco-russe en math\'{e}matiques'' du CNRS }

\address{
Manfred Lehn\\
Fachbereich Physik, Mathematik und Informatik\\
Johannes Gutenberg--Uni\-ver\-si\-t\"{a}t Mainz\\
D-55099 Mainz, Germany}
\email{lehn@mathematik.uni-mainz.de}

\address{
Christoph Sorger\\
Laboratoire de Math\'{e}matiques Jean Leray (UMR 6629 du CNRS)\\
Universit\'{e} de Nantes\\
2, Rue de la Houssini\`{e}re\\
BP 92208\\
F-44322 Nantes Cedex 03, France}
\email{christoph.sorger@univ-nantes.fr}

\subjclass{Primary 14J60; Secondary 14D20, 14J28, 32J27}

\maketitle

%%%%%%%%%%%%%%%%%%%%%%%%%%%%%%%%%%%%%%%%%%%%%%%%%%%%%%%%%%%%%%%%%%%%%%%%%%%%%
%%%%%%%%%%%%%%%%%%%%%%%%%%%%%%%%%%%%%%%%%%%%%%%%%%%%%%%%%%%%%%%%%%%%%%%%%%%%%

\section{Introduction}

%%%%%%%%%%%%%%%%%%%%%%%%%%%%%%%%%%%%%%%%%%%%%%%%%%%%%%%%%%%%%%%%%%%%%%%%%%%%%
%%%%%%%%%%%%%%%%%%%%%%%%%%%%%%%%%%%%%%%%%%%%%%%%%%%%%%%%%%%%%%%%%%%%%%%%%%%%%
How to construct irreducible holomorphic symplectic manifolds?
Except for the variety of Beauville and Donagi \cite{BeauvilleDonagi}
all known examples arise from moduli spaces of semistable sheaves
on a K3 or abelian surface.

For every element $v$ in the Mukai lattice $H^{\even}(X,\IZ)$ of a polarised
K3 or abelian surface $(X,H)$ there is an associated moduli space $M_v$ that
parametrises polystable sheaves $E$ with Mukai vector $v=v(E):=ch(E)\sqrt{\td(X)}$.
If $H$ and $v$ are chosen to the effect that no strictly semistable sheaves exist,
i.e.\ every semistable sheaf is automatically stable, then $M_v$ is
a projective holomorphically symplectic manifold due to Mukai \cite{Mukai}.

In the opposite case, $M_v$ is singular and one may ask whether $M_v$ at
least admits a projective symplectic resolution. This question has
been raised and successfully answered in two cases by O'Grady \cite{OG1,OG2},
leading to two new deformation classes of irreducible holomorphic symplectic
manifolds.

In this paper we give a complete answer to O'Grady's question for general
ample divisors $H$ and moduli spaces whose expected dimension $2+\langle v,v\rangle$ is $\geq 4$.
The answer depends essentially only on the divisibility of the Mukai
vector $v\in H^{\even}(X,\IZ)$ and the dimension of the moduli space.
We may write $v=mv_0$ with a primitive Mukai vector $v_0=(r,c,a)$ and
a multiplicity $m\in\IN$. Suppose for simplicity that $r>0$, and let
$H$ denote a $v$--general ample divisor. Then every semistable sheaf $E$ with
Mukai vector $v(E)=v_0$ is stable, and a necessary and sufficient
condition for the existence of $E$ is that $c\in\NS(X)$ and that
$\langle v_0,v_0\rangle \geq -2$. There are five principal cases to distinguish:

1) If $\langle v_0,v_0\rangle=-2$, then Mukai has shown that $M_{v_0}$
consists of a single point $[E_0]$ only. As the expected dimension of the
moduli space $M_{mv_0}$ is negative for $m>0$ there are no stable sheaves in
this case, and it follows by induction that any semistable sheaf must be of
the form $E_0^{\oplus m}$. Hence $M_v$ is a single point as well.

2) If $\langle v_0,v_0\rangle=0$, the moduli space $M_{v_0}$ is
again a K3 surface or an abelian surface if $X$ is K3 or abelian due
to beautiful results of Mukai. It turns out that any semistable
sheaf $E$ with $v(E)=mv_0$ is $S$-equivalent to a direct sum
$E=E_1\oplus \ldots\oplus E_m$ with stable sheaves $[E_i]\in
M_{v_0}$. It follows that $M_v=S^m(M_{v_0})$. Thus the moduli spaces
are singular in codimension 2, but admit symplectic resolutions in
terms of the Hilbert scheme $\Hilb^m(M_{v_0})\to M_v$.

3) Assume now that $\langle v_0,v_0\rangle \geq 2$. Due to the
combined efforts of many authors, with important steps taken by
Mukai, Huybrechts, O'Grady and Yoshioka, one finally has the
following result \cite{Yoshioka}:  $M_{v_0}$ is a smooth symplectic
variety that is deformation equivalent to $\Hilb^{1+\tfrac12\langle
v_0,v_0\rangle}(X)$, if $X$ is a K3-surface, and to $\Pic_0(X)\times
\Hilb^{\tfrac12\langle v_0,v_0\rangle}$, if $X$ is an abelian
surface.

Assume in addition that $m\geq 2$. The main result of this article implies
that one has to further distinguish the following two cases:

4) Let $\langle v_0,v_0\rangle =2$ and $m=2$. The moduli spaces $M_{K3}(2;0,4)$
and $M_{Ab}(2;0,2)$ studied by O'Grady \cite{OG1,OG2} and Rapagnetta \cite{Rapagnetta} fall into this class.
The moduli space $M_v$ has dimension $10$, its singular locus has codimension 2
and is in fact isomorphic to $S^2M_{v_0}$. As shown in \cite{LS4}, the
symplectic desingularisations constructed by O'Grady exist for all Mukai vectors
in this class and can be obtained by blowing-up the reduced singular locus.

5) In all other cases our main result states:\\

\noindent
{\bf Theorem A} --- {\sl If either $m\geq 2$ and $\langle v_0,v_0\rangle >2$
or $m>2$ and $\langle v_0,v_0\rangle\geq 2$, then $M_{mv_0}$ is a
locally factorial singular symplectic variety.}\\

As an immediate application one obtains:\\

\noindent
{\bf Theorem B} --- {\sl Under the hypotheses of Theorem A, $M_{mv_0}$ does not
admit a proper symplectic resolution.}\\

Under some technical hypotheses theorems A and B hold as well for semistable
torsion sheaves (see the main text).
Partial results for Theorem B in the case $m=2$ have been obtained
previously by two of us \cite{KL}, and, independently and with different
methods, by Kiem and Choy \cite{Kiem,KiemChoy}.

We note that our approach is rather general; our main technical
result, Proposition \ref{prop:sympredfactorial}, is essentially a
linear-algebraic fact. Therefore, we expect that results similar to
Theorem A and B might hold in other situations with similar geometry
-- in particular, for the moduli spaces of flat connections on an
algebraic curve. In fact, Proposition \ref{prop:sympredfactorial} is
a statement about quiver varieties of H.~Nakajima \cite{Nakajima},
although the numerical data corresponding to our quivers are
specifically excluded from  consideration in \cite{Nakajima} (which
is not surprising as one of the results of \cite{Nakajima} is that
any quiver variety considered there {\sl does} admit a symplectic
resolution). Thus our approach and our Proposition
\ref{prop:sympredfactorial} might be used wherever one finds quiver
varieties of the
same type.\\

Acknowledgements: The moment when Grothendieck's theorem on
factoriality was shown to us by Duco van Straten turned around our
approach to the problem. We thank him as well as Stefan Bauer,
Daniel Huybrechts and Joseph Le Potier for many helpful discussions.

Many authors have worked on moduli of sheaves on K3 and abelian surfaces ever since
the seminal work of Mukai. The most general results for our purposes have been
obtained by Yoshioka \cite{Yoshioka}. We refer to the textbook \cite{HL} and
Yoshioka's paper for further references and general information on semistable
sheaves and their moduli spaces.

%%%%%%%%%%%%%%%%%%%%%%%%%%%%%%%%%%%%%%%%%%%%%%%%%%%%%%%%%%%%%%%%%%%%%%%%%%%%%%%%%%%
%%%%%%%%%%%%%%%%%%%%%%%%%%%%%%%%%%%%%%%%%%%%%%%%%%%%%%%%%%%%%%%%%%%%%%%%%%%%%%%%%%%
\section{Notation and conventions, plan of the paper}
%%%%%%%%%%%%%%%%%%%%%%%%%%%%%%%%%%%%%%%%%%%%%%%%%%%%%%%%%%%%%%%%%%%%%%%%%%%%%%%%%%%
%%%%%%%%%%%%%%%%%%%%%%%%%%%%%%%%%%%%%%%%%%%%%%%%%%%%%%%%%%%%%%%%%%%%%%%%%%%%%%%%%%%

\subsection{The underlying surface.}
Throughout this article $X$ will denote a complex projective K3 or abelian
surface with a fixed symplectic structure, i.e.\ an isomorphism
$H^2(X,\ko_X)\isom\IC$, and a fixed ample divisor $H$.

The even integral cohomology $H^{\even}(X,\IZ)$ is equipped with a pairing
$$\langle v,w\rangle:=-\int_X vw^\dual,$$
where $w^\dual=(-1)^iw$ for $w\in H^{2i}(X,\IZ)$. Following Mukai we associate
to each coherent sheaf $E$ its Mukai vector
$$v(E):=\ch(E)\sqrt{\td(X)}\in H^{\even}(X,\IZ).$$
The Hilbert polynomial of $E$ with respect to an ample divisor $H$
can be expressed in terms of its Mukai vector as follows:
$$\chi(E\tensor\ko_X(mH))=-\langle v(E), v(\ko_X(-mH))\rangle=:P_v(m).$$

\subsection{Semistable sheaves.}
Stability or semistability of a coherent sheaf is defined with respect
to a fixed ample divisor $H$. We let $M_v$ denote the moduli space of
semistable sheaves with Mukai vector $v$. Closed points of $M_v$ are in natural
bijection with polystable sheaves $E$. Points corresponding to stable sheaves form
a -- possibly empty -- open subset $M_v^s\subset M_v$.

Semistable sheaves may have two--, one-- or zero--dimensional support.
Stability in the first case was defined by Maruyama and Gieseker, the
generalisation to pure sheaves of arbitrary dimension is due to Simpson.
In the rest of the paper we exclude once for all the case of zero-dimensional
sheaves as being well-known: if the Mukai vector is $v=(0,0,a)$ then
$M_v\isom S^aX$, the symmetric product of $X$, and the Hilbert-Chow morphism
$\Hilb^a(X)\to M_v$ provides a symplectic resolution.

\subsection{General assumptions.}\label{subsec:GeneralAssumptions}
An element $v_0\in H^{\even}(X,\IZ)$ is primitive if it is not an integral
multiple of another lattice element. Given a non-trivial element
$v\in H^{\even}(X,\IZ)$ we may always decompose it as $v=mv_0$ with a primitive
element $v_0$ and a multiplicity $m\in\IN$. Throughout this article we assume
that $v_0=(r_0,c_0,a_0)$ has the following properties:
\begin{center}
$(*)$\quad\quad\quad$\left\{
\begin{array}{l}
\text{Either $r_0>0$ and $c_0\in\NS(X)$,}\\
\text{\quad\quad or $r_0=0$, $c_0\in\NS(X)$ is effective, and $a_0\neq 0$;}\\
\text{$\langle v_0,v_0\rangle \geq 2$.}
\end{array}
\right.$
\end{center}

The results of this paper suggest to distinguish systematically between
the following three cases for a Mukai vector $v$ satisfying assumptions $(*)$:
\begin{itemize}
\item[(A)] $m=1$.
\item[(B)] $m=2$ and $\langle v_0,v_0\rangle=2$.
\item[(C)] $m\geq 3$, or $m=2$ and $\langle v_0,v_0\rangle \geq 4$.
\end{itemize}

%%%%%%%%%%%%%%%%%%%%%%%%%%%%%%%%%%%%%%%%%%%%%%%%%%%%%%%%%%%%%%%%%%%%%%%%%%%%%%%%%%%%%%
\subsection{General ample divisors.}\label{subsec:GeneralAmpleDivisors}
The significance of $(*)$ lies in the
fact that one has the notion of a $v$--general ample divisor $H$: there is a
systems of
hyperplanes in the ample cone of $X$, called $v$--walls, that is countable
but locally finite for torsion free sheaves (\cite{HL}, ch.\ 4C) and finite for
torsion sheaves (\cite{Yoshioka}, sec.\ 1.4.) with the following
property: if $H$ is $v$--general, i.e.\ if $H$ is not contained in any
$v$--wall, then for every direct summand $E'$ of a polystable sheaf $E$
with $v(E)=v$ one has $v(E')\in \IQ v(E)$.

Let $H$ be a $v_0$--general ample divisor and consider the following
assertions:
$$
\begin{array}{rl}
  (**) & \text{$M_{v_0}$ is non-empty.} \\
  (***) & \text{$M_{v_0}$ is irreducible.}
\end{array}
$$
Yoshioka shows in \cite{Yoshioka}, Thm 0.1 and Thm 8.1, that $(*)$
implies $(**)$ and $(***)$ except when $X$ is a K3 surface, $r_0=0$
and $c_0$ is not ample. Moreover he has communicated to us an
unpublished note that fills this gap, so that $(**)$ and $(***)$ are
consequences of $(*)$. An essential technique in Yoshioka's work is
the deformation of the underlying surface; the arguments are rather
involved. For the irreducibility part $(***)$ we give a new and
direct proof, based on an old and beautiful idea of Mukai, see
Theorem \ref{thm:irreducibility}.

%%%%%%%%%%%%%%%%%%%%%%%%%%%%%%%%%%%%%%%%%%%%%%%%%%%%%%%%%%%%%%%%%%%%%%%%%%%%%%%%%%%%%%
\subsection{Elements of the construction of moduli spaces.}
\label{subsec:constructionofmoduli}
We need to recall some basic elements of the construction machinery of moduli
spaces of sheaves following the approach of Simpson \cite{Simpson} (see also
\cite{HL}, ch.\ 4). Let $v$ be a Mukai vector satisfying $(*)$ and let
$P_v$ denote the corresponding Hilbert polynomial. Choose a sufficiently large
integer $k=k(v)$ and put $N=P_v(k)$, $\kh:=\ko_X(-kH)^{\oplus N}$. Then there
is a closed subscheme $R\subset\Quot_{X,H}(\kh,P)$ with the following property:
a closed point
$$[q:\kh\to E]\in R$$
is stable or semistable with respect to the canonical $\LiePGl(N)$--action
and the corresponding linearisation of the determinant line bundle on $R$
if and only if $q$ induces an isomorphism $\IC^{N}\to H^0(X,E(kH))$ and
if $E$ is stable or semistable, respectively.
Let $R^{s}\subset R^{ss}\subset R$ denote the open subsets of stable and
semistable points, respectively. Then
$$R^{ss}\GIT\LiePGl(N))\isom M_v\quad\quad\text{and}\quad\quad
R^{s}\GIT\LiePGl(N)\isom M^s_v.$$
Let $\pi: R^{ss}\to M_v$ denote the quotient map.
The orbit of a point $[q:\kh\to E]$ is closed in $R^{ss}$
if and only if $E$ is polystable. In that case, the stabiliser subgroup of $[q]$ in
$\LiePGl(N)$ is canonically isomorphic to $\PAut(E)=\Aut(E)/\IC^*$.
Moreover, by Luna's slice theorem there is a $\PAut(E)$--invariant subscheme
$S\subset R^{ss}$, containing $[q]$, such that the canonical morphisms
$$(\LiePGl(N)\times S)\GIT \PAut(E)\to R^{ss}\quad\quad\text{and}\quad\quad
S\GIT\PAut(E)\to M$$
are \'etale. The Zariski tangent space $T_{[q]}S$ is isomorphic to $\Ext^1(E,E)$.

%%%%%%%%%%%%%%%%%%%%%%%%%%%%%%%%%%%%%%%%%%%%%%%%%%%%%%%%%%%%%%%%%%%%%%%%%%%%%%%%%%%%%
\subsection{Local description.}\label{subsec:LocalDescription}
The completion of the local ring $\ko_{S,[q]}$ has the following
deformation theoretic description:

Let $\IC[\Ext^1(E,E)]$ denote the ring of polynomial functions on
$\Ext^1(E,E)$ and let $A:=\IC[\Ext^1(E,E)]^{\wedge}$ denote its
completion at the maximal ideal $\mathfrak{M}$ of functions
vanishing at $0$. There is a trace map $\tr:\Ext^2(E,E)\to
H^2(\ko_X)$. We denote its kernel by $\Ext^2(E,E)_0$. The
automorphism group $\Aut(E)$ naturally acts on $\Ext^1(E,E)$ and
$\Ext^2(E,E)_0$ by conjugation. Since the scalar multiples of the
identity act trivially we actually have an action of the projective
automorphism group $\PAut(E)=\Aut(E)/\IC^{*}$. There is a linear map
$$\kappa:\Ext^2(E,E)_0^*\lra\IC[\Ext^1(E,E)]^{\wedge},$$
the so-called Kuranishi map, with the following properties:
\begin{enumerate}
\item $\kappa$ is $\PAut(E)$--equivariant.
\item Let $I$ be the ideal generated by the image of $\kappa$. Then
there are isomorphisms of complete rings
$$\widehat\ko_{S,[q]}\isom A/I\quad
\text{and}\quad
\widehat\ko_{M_v,[E]}\isom (A/I)^{\PAut(E)}.$$
\item\label{it:property3}
For every linear form $\varphi\in\Ext^2(E,E)_0^*$ one has,
for $e\in \Ext^1(E,E)$,
$$\kappa(\varphi)(e)=\tfrac12\varphi(e\cup e)+\text{higher order terms in }e.$$
\end{enumerate}

%%%%%%%%%%%%%%%%%%%%%%%%%%%%%%%%%%%%%%%%%%%%%%%%%%%%%%%%%%%%%%%%%%%%%%%%%%%%%%%%%%%%%%%%%
\subsection{Passage to the normal cone.}\label{subsec:PassageToNormalCone}
Let $J\subset \IC[\Ext^1(E,E)]$ denote the ideal generated by the image of
the quadratic part of $\kappa$:
$$\kappa_2:\Ext^2(E,E)_0^*\lra S^2\Ext^1(E,E)^*,\quad \varphi\mapsto(e\mapsto\tfrac12\varphi(e\cup e))$$
Then $J$ is the ideal of the null-fibre $F=\mu^{-1}(0)$ of the morphism
$$\mu:\Ext^1(E,E)\lra\Ext^2(E,E)_0,\quad \mu(e)=\tfrac12(e\cup e).$$
The ideals $I\subset \IC[\Ext^1(E,E)]^{\wedge}$ and $J\subset \IC[\Ext^1(E,E)]$
are related as follows. The graded ring $\gr A$ associated to the $\gothm$-adic
filtration on $A=\IC[\Ext^1(E,E)]^{\wedge}$ is canonically isomorphic to
$\IC[\Ext^1(E,E)]$. For any ideal $\gotha\subset A$ let
$\ini(\gotha)\subset \gr A$ denote the ideal generated by the leading terms
$\ini(f)$ of all elements $f\in\gotha$. Then property \ref{it:property3} of the
Kuranishi--map says that
$$J\subset \ini(I).$$
Hence there is the following chain of inequalities:
\begin{equation}\label{eq:dimestimate}
\begin{array}{rcl}
\dim(F)&=&\dim\;(\gr A)/J\\
&\geq &\dim\;(\gr A)/\ini(I)\;=\;\dim\;\gr (A/I)\;=\;\dim(A/I)\\
&\geq& \dim \Ext^1(E,E)-\dim \Ext^2(E,E)_0,
\end{array}
\end{equation}
where the last inequality comes from the fact that $A$ is regular
of dimension $=\dim \Ext^1(E,E)$ and $I$ is generated by
$\dim \Ext^2(E,E)_0$ elements.

We need to describe $\mu$ in greater detail; the resulting description is similar
to Nakajima's construction of the so-called quiver varieties \cite{Nakajima}. Write
\begin{equation}
E=\bigoplus_{i=1}^s W_i\tensor E_i
\end{equation}
with pairwise non-isomorphic stable sheaves $E_i$ and vector spaces $W_i$ of dimension
$n_i$. Let $W_{ij}:=\Hom(W_i,W_j)$ and
$V_{ij}:=\Ext^1(E_i,E_j)$. Then
$$\End(E)=\bigoplus_iW_{ii},\quad \Ext^1(E,E)=\bigoplus_{i,j}W_{ij}\tensor V_{ij},
\quad \Ext^2(E,E)=\bigoplus_iW_{ii}.$$
The automorphism group
$$\Aut(E)=\prod_i\Aut(W_i)\isom\prod_i\LieGl(n_i)=:G(n)$$
acts on $\Ext^1(E,E)$ by conjugation on the first factor in each direct
summand. By Serre-Duality, the pairing
$$V_{ij}\tensor V_{ji}\to \IC,\quad e\tensor e'\mapsto \tr(e'\cup e)$$
is non-degenerate and antisymmetric. This yields a symplectic form $\omega$ on
$\Ext^1(E,E)$ such that $W_{ij}\tensor V_{ij}$ and $W_{ab}\tensor V_{ab}$
are perpendicular, unless $i=b$ and $j=a$, in which case
$$\omega:(W_{ij}\tensor V_{ij})\tensor (W_{ji}\tensor V_{ji})\lra \IC,
\quad \omega(A\tensor e, A'\tensor e')= \tr(A'A)\tr(e'\cup e).$$
Moreover, the quadratic map $\mu:\Ext^1(E,E)\to \Ext^2(E,E)_0$ is given by
\begin{equation}\label{eq:muexplizit}
\mu\left(\sum_{ij} \sum_k A_{ij}^k\tensor e_{ij}^k\right)
=\sum_{ij}\sum_{k,\ell} A_{ij}^kA_{ji}^\ell \tr(e_{ij}^ke_{ji}^\ell).
\end{equation}

%%%%%%%%%%%%%%%%%%%%%%%%%%%%%%%%%%%%%%%%%%%%%%%%%%%%%%%%%%%%%%%%%%%%%%%%%%%%%%%%%%%%%%%%
\subsection{Strategy}\label{subsec:strategy}
In general we do not know how to compute the Kuranishi map explicitly. However,
the explicit description of the quadratic part $\mu$ given above allows
for a detailed study of the fibre $F:=\mu^{-1}(0)\subset \Ext^1(E,E)$.
The passage from $\kappa$ to $\mu$ corresponds to the passage from the local ring
$\widehat\ko_{S,[q]}$ to the coordinate ring $\ko_F$ of its tangent cone.

In section \ref{sec:symplectic} we show that under certain
hypotheses the fibre $F$ is an irreducible normal complete
intersection which is, in case (C), regular in codimension $\leq 3$
and state consequences for the local rings $\ko_{R^{ss},[q]}$ of
points $[q]$ in closed orbits of $R^{ss}$.

Section \ref{sec:irreducibility} contains a basic irreducibility
result for moduli spaces of sheaves on a K3 or abelian surface.

In section \ref{sec:Misfactorial} it is proved --- under the
hypothesis that the ample divisor is $v$--general --- that the
moduli space $M$ is a non-empty irreducible normal variety of
expected dimension, and that it is locally factorial in case (C). As
an application we show in section \ref{sec:SymplecticResolution}
that in case (C) the moduli space does not admit a symplectic
resolution.

%%%%%%%%%%%%%%%%%%%%%%%%%%%%%%%%%%%%%%%%%%%%%%%%%%%%%%%%%%%%%%%%%%%%%%%%%%%%%
%%%%%%%%%%%%%%%%%%%%%%%%%%%%%%%%%%%%%%%%%%%%%%%%%%%%%%%%%%%%%%%%%%%%%%%%%%%%%

\section{Symplectic reduction}\label{sec:symplectic}

%%%%%%%%%%%%%%%%%%%%%%%%%%%%%%%%%%%%%%%%%%%%%%%%%%%%%%%%%%%%%%%%%%%%%%%%%%%%%
%%%%%%%%%%%%%%%%%%%%%%%%%%%%%%%%%%%%%%%%%%%%%%%%%%%%%%%%%%%%%%%%%%%%%%%%%%%%%

\subsection{The symplectic momentum map.}
Let $U$ be a smooth affine algebraic variety over $\IC$ endowed with
a symplectic form $\omega$. Let $G$ be a reductive group that acts
on $U$ preserving $\omega$. This action induces an infinitesimal
action of the Lie algebra $\gothg$ of $G$, i.e.\ a homomorphism of
Lie algebras $\gothg\to \Gamma(U,T_U)$. We denote the vector field
corresponding to $A\in\gothg$ at $x\in U$ by $A_x$. A momentum map
for the action is a $G$-equivariant morphism $\mu:U\to\gothg^*$ with
the property that $d\mu_x(\xi)(A)=\omega(\xi,A_x)$ for all $x\in U$
and $\xi\in T_xU$. If a momentum map exists, it is unique up to an
additive constant in $(\gothg^*)^G$.

Let $\mu:U\to \gothg^*$ be a momentum map with null-fibre $F:=\mu^{-1}(0)$.

\begin{lemma}\label{lem:muimage}---
Let $x\in F$ be a point with stabiliser subgroup $H\subset G$.
Then the image of $d\mu_x:T_xU\to \gothg^*$ is $(\gothg/\gothh)^*=\gothh^\perp$,
where $\gothh\subset\gothg$ denotes the Lie algebra of $H$. In particular,
if $H$ is finite then $d\mu_x$ has maximal rank and $F$ is regular at $x$ of
dimension $\dim(U)-\dim(G)$.
\end{lemma}

\begin{proof} The image $d\mu_x$ annihilates $A\in\gothg$ if and only if
$\omega(\xi,A_x)=0$ for all $\xi\in T_xU$, i.e.\ if $A_x$ is perpendicular
to $T_xU$ with respect to $\omega$. As $\omega$ is non-degenerate, this
is equivalent to saying that $A_x$ vanishes, hence is a tangent vector
to the stabiliser subgroup $H$.
\end{proof}

\begin{lemma}\label{lem:boundssuffice}--- Let $\mu:U\to \gothg^*$ be a
momentum map with null-fibre $F$. Let $Z\subset F$ be the closed subset of
points with non-finite stabiliser group. Let $d=\dim U-\dim\gothg$.
\begin{enumerate}
\item If $\dim(Z)\leq d-1$, then $F$ is a reduced complete intersection
of dimension $d$.
\item If $\dim(Z)\leq d-2$, then $F$ is normal.
\end{enumerate}
\end{lemma}

\begin{proof} Every irreducible component of $F$ must have dimension
$\geq d$ since $F$ is cut out by $\dim\gothg$ equations. By Lemma
\ref{lem:muimage}, $F$ has dimension $d$ in each point $x\in
F\setminus Z$. If $\dim Z<d$, then $F\setminus Z$ is dense in $F$
and every irreducible component has precisely dimension $d$. Hence
$F$ is a complete intersection and in particular Cohen-Macaulay
(\cite{AK}, Cor.\ III 4.5). Since $F\setminus Z$ is smooth, $F$ in
addition satisfies condition $R_0$ and is therefore reduced
(\cite{AK}, Prop.\ VII 2.2). If in addition $\dim(Z)\leq d-2$, then
$F$ is regular in codimension 1 and normal by Serre's criterion
(\cite{AK}, Cor.\ VII 2.13).
\end{proof}

%%%%%%%%%%%%%%%%%%%%%%%%%%%%%%%%%%%%%%%%%%%%%%%%%%%%%%%%%%%%%%%%%%%%%%%%%%%%%%%%%%
\subsection{The key estimate.}
We want to apply the lemma to the following particular situation, that arises
in the study of local rings of the moduli space of sheaves.

{\sl Set-up:} Let $W_1,\ldots,W_s$ be a sequence of vector spaces,
$s\geq 1$. The dimensions $n_i=\dim(W_i)$ form the components of a
vector $n\in\IN^{s}_0$. Furthermore, let $W_{ij}=\Hom(W_i,W_j)$.
There is a natural symmetric pairing
$$W_{ij}\tensor W_{ji}\to \IC,\quad\quad(A,B)\mapsto \tr(BA).$$
Moreover, let $V_{ij}$, $1\leq i,j\leq s$, be vector spaces, equipped with
non-degenerate pairings
$$\omega_{ij}:V_{ij}\tensor V_{ji}\to \IC,$$
that are skew-symmetric in the sense that $\omega_{ij}(e,e')=-\omega_{ji}(e',e)$.
Then the vector space $U(n):=\bigoplus_{i,j} W_{ij}\tensor V_{ij}$ carries
a natural symplectic form $\omega$ with the property that $W_{ij}\tensor V_{ij}$
is perpendicular to all $W_{ab}\tensor V_{ab}$, $(a,b)\neq (j,i)$ and
$$\omega:\Big(W_{ij}\tensor V_{ij}\Big)\tensor \Big(W_{ji}\tensor V_{ji}\Big)
\lra \IC,\quad
(A\tensor e)\tensor (A'\tensor e')\mapsto \tr(A'A)\omega_{ij}(e,e').
$$
In the following arguments the vector spaces $V_{ij}$ are fixed and chosen
once for all, whereas the sequence of vector spaces $W_i$ can be replaced
by appropriate subspaces etc. We will argue by induction
over the dimension vector $n$ as an element in the monoid $\IN_0^s$.
Most objects defined below will therefore be indexed by $n$, like the space
$U(n)$ above, even if this is not quite accurate as they really depend on the
spaces $W_i$.

The group $G(n)=\prod_i\Aut(W_i)$ acts on $U(n)$ by conjugation on
the first factors in the decomposition. The subgroup of scalars
$\IC^*\subset G(n)$ acts trivially. Let $PG(n):=G(n)/\IC^*$.
The action of $PG(n)$ on $U(n)$ preserves the symplectic structure.
The moment map for the action is
\begin{eqnarray*}
\mu(n):U(n)&\lra& \gothp\gothg(n)^*\simeq\Kern\Big(\bigoplus_i\Liegl(n_i)\xra{\;\;\tr\;\;}\IC\Big),\\
\sum_{i,j,k}A_{ij}^{(k)}v_{ij}^{(k)}
&\mapsto&
\sum_{k,\ell}\sum_{i,j}
A_{ij}^{(k)}A_{ji}^{(\ell)}\tr(v_{ij}^{(k)}\cup v_{ji}^{(\ell)})
\end{eqnarray*}

Let $F(n):=\mu(n)^{-1}(0)\subset U(n)$ denote the null-fibre of the moment map.
The structure of $F(n)$ depends only on $n$ and the dimensions $d_{ij}:=\dim(V_{ij})$.
Let $D$ denote the matrix $(d_{ij})$ and let $a:=\min\{d_{ij}-2\delta_{ij}\}$.

\begin{proposition}\label{prop:sympredfactorial}---
Assume that $a\geq 2$. Then $F(n)$ is an irreducible normal complete
intersection of dimension $d:=n\tra (D-I) n+1$. Moreover, $F(n)$ is
regular in codimension $\leq 3$ with the possible exception of the
two cases
\begin{enumerate}
\item $n=(1,1)$, $d_{12}=2$, and
\item $n=(2)$, $d_{11}=4$.
\end{enumerate}
\end{proposition}

\begin{proof} 1. Since $\dim(U(n))=\sum_{i,j}n_in_j d_{ij}$ and since the
range of $\mu$ has dimension $\sum_{i}n_in_i-1$, the expected dimension
of $F(n)$ is
$$d=\sum_{i,j}n_in_j d_{ij}-\sum_{i}n_in_i+1=n\tra (D-I)n+1.$$
Also, $F(n)$ is a cone and hence connected. By Lemma \ref{lem:boundssuffice},
it suffices to show that the locus $Z$ of
points in $F(n)$ with non-trivial stabiliser in $PG(n)$ has dimension
$\leq d-4$ in general and $\leq d-3$ in the two exceptional cases. This
will be done by induction on the dimension vector $n\in\IN^{s}_0$.

The induction starts with $n=(0,\ldots,1,\ldots,0)$, in which case the
statement is trivial. So let $n\in\IN_0^s$ be an arbitrary element and
assume that the proposition holds for all $n'\in\IN^s_0$ such
that $0<\sum_i{n_i'}<\sum_i{n_i}$.

2. We can analyse $Z$ as follows: Let $g\in G(n)$, $g\notin\IC^*$,
and consider the corresponding fixed point locus $F(n)^g$. The image
$G(n) F(n)^g$ of the morphism $\varphi:G(n)\times F(n)^g\to F(n)$,
$(g',x)\mapsto g'x$, consists of all points $y\in F$ whose
stabiliser subgroup $G(n)_y$ contains an element conjugate to $g$.
Suppose that $H\subset G(n)$ is a subgroup that stabilises the fixed
point set $F(n)^g$.  Then we can bound the dimension of the fibres
of $\varphi$ by $\dim(H)$. It follows that $\dim(G(n)
F(n)^g)\leq\dim F(n)^g+\dim G(n)-\dim H$. In the following we will
describe a finite set of elements $g$ such that $Z$ is covered by
the corresponding sets $G(n)F(n)^{g}$ and such that for each $g$ one
has $\dim(G(n)F(n)^g)\leq d-3$ or $\leq d-4$. This gives the desired
bound for $\dim(Z)$.

3. Let $g=(g_1,\ldots,g_s)\in G(n)$, $g\notin\IC^*$. We distinguish three cases:

3.1. {\sl Case: $g$ is semisimple.} For each $\lambda\in\IC$ consider the
eigenspaces $W_i(\lambda)\subset W_i$ of $g_i$, and let $n_i(\lambda)=
\dim(W_i(\lambda))$, $n(\lambda)=(n_i(\lambda))_i$. Then
$n=\sum_\lambda n(\lambda)$. There is a decomposition
$$U(n)^g=\bigoplus_\lambda U(n(\lambda)),\quad
U(n(\lambda))=\bigoplus_{ij}\Hom(W_i(\lambda),W_j(\lambda))\tensor V_{ij}.$$
Moreover, the restriction of the momentum map to the fixed point locus splits
into a product of momentum maps for each $U(n(\lambda))$:
$$\mu(n)|_{U(n)^g}=\prod_\lambda \mu(n(\lambda)),\quad
\mu(n(\lambda)):U(n(\lambda))\lra \gothp\gothg(n(\lambda))^*.$$ It
follows that
$$F(n)^g=\prod_\lambda F(n(\lambda))\quad\text{with }
F(n(\lambda))=\mu(n(\lambda))^{-1}(0).$$
By induction, we have
$$\dim(F(n)^g)=\sum_{\lambda} \dim(F(n(\lambda))=\sum_{\lambda}{}'
\Big(n(\lambda)\tra (D-I)n(\lambda)+1\Big),$$
where $\sum{}'$ indicates that we only sum over all $\lambda$ with $n(\lambda)\neq 0$.

Next, $U(n)^g$ is stabilised by $H=\prod_\lambda G(n(\lambda))$, a subgroup
in $G(n)$ of codimension $n\tra n-\sum_\lambda n(\lambda)\tra n(\lambda)$. We
obtain the following upper bound for the dimension of $G(n)F(n)^g$:
$$\dim(G(n)F(n)^g)\leq \sum_\lambda{}'\Big(n(\lambda)\tra(D-2I) n(\lambda)+1\Big)+n\tra n.$$
Note that $\nu:=|\{\lambda\;|\;n(\lambda)\neq 0\}|\geq 2$, since $g\notin\IC^*$.
The difference of $\dim(G(n)F(n)^g)$
to the expected dimension of $F(n)$ is therefore bounded below by
\begin{eqnarray*}
\Delta&:=&\big( n\tra (D-2I)n+1\big)
-\sum_\lambda{}'\big( n(\lambda)\tra(D-2I)n(\lambda)+1\big)\\
&=&\sum_{\lambda\neq \mu} n(\lambda)\tra(D-2I)n(\mu) - (\nu-1)\\
&\geq&2\nu(\nu-1)-(\nu-1)=(2\nu-1)(\nu-1)\geq 3.
\end{eqnarray*}
Clearly, $\Delta\geq 4$ for $\nu\geq 3$. Assume that $\nu=2$, say with the
distinct eigenvalues $\lambda$ and $\lambda'$. Then
$$\Delta\geq 2\sum_{i,j} n_i(\lambda')n_j(\lambda'')(d_{ij}-2\delta_{ij})-1
\geq 2a \sum_i n_i(\lambda')\sum_i n_i(\lambda)-1.$$
Thus $\Delta=3$ implies $a=2$ and $\sum_in_i(\lambda)=1=\sum_in_i(\lambda')$.
Hence there are only the following exceptional cases:
\begin{enumerate}
\item $s=1$, $n=2$, $d_{11}=2+2\delta_{11}=4$, or
\item $s=2$, $n=(1,1)$ and $d_{12}=d_{21}=2$.
\end{enumerate}

If a point $f \in F(n)$ is fixed by a semisimple element, it is also
fixed by a whole subtorus $T \subset G(n)$. Up to a conjugation,
there is only a finite number of such subtori $T_i \subset
G(n)$. Choosing an element $g_i \in T_i$ in each of these subtori,
we see that the union of all sets $G(n)F(n)^g$, $g$ semisimple, is
covered by the finite union of all sets $G(n)F(n)^{g_i}$.

3.2. {\sl Case: $g$ is unipotent.} We may write $g=1+h$, with a non-zero
nilpotent element $h=(h_1,\ldots,h_s)\in\bigoplus_i \End(W_i)$.
Let $K_i^{(\ell)}:=\ker(h_i^\ell)\subset W_i$ and $m_i(\ell):=\dim K_i^{(\ell)}$
for all $\ell\in\IN_0$. There is a filtration
$$0=K_i^{(0)}\subset K_i^{(1)}\subset\ldots =W_i.$$
For each level $\ell>0$ we choose a graded complement $W_i^{(\ell)}$ to
$h K_i^{(\ell+1)}+K_i^{(\ell-1)}$ in $K^{(\ell)}$ and let
$n_i^{(\ell)}=\dim W_i^{(\ell)}$. (We note that this is an instance
of the so-called Jacobson-Morozov-Deligne filtration associated to a
nilpotent element, see \cite[1.6]{WeilII}; the spaces $W_i^{(\ell)}$
are the primitive subspaces with respect to an $sl_2$-triple
containing $h$.)

Suppose that $A=(A_{ij})\in \bigoplus_{ij}\Hom(W_i,W_j)$ commutes with $h$.
Then $A_{ij}$ is completely determined by its value on the spaces $W_i^{(\ell)}$,
$\ell\in\IN$, and conversely, any value of $A_{ij}:W_i^{(\ell)}\to K_j^{(\ell)}$
can be prescribed. The composition with the canonical projection
$K_j^{(\ell)}\to W_j^{(\ell)}$ defines a homomorphism
$A_{ij}^{(\ell)}:W_i^{(\ell)}\lra W_j^{(\ell)}$,
and the map
$$\Phi: \left(\bigoplus_{ij}\Hom(W_i,W_j)\right)^g\lra \bigoplus_{\ell}
\left(\bigoplus_{ij}\Hom(W_i^{(\ell)},W_j^{(\ell)})\right),\quad
(A_{ij})\mapsto (A_{ij}^{(\ell)}),$$
is a ring homomorphism. Let
$$\Phi_V:U(n)^g=\left(\bigoplus_{i,j}W_{ij}\tensor V_{ij}\right)^g\lra
\bigoplus_{\ell} U(n^{(\ell)})=\bigoplus_{\ell}\bigoplus_{i,j}\Hom(W_i^{(\ell)},
W_j^{(\ell)})\tensor V_{ij}$$
be analogously defined. Then $\Phi_V(F(n)^g)\subset\prod_\ell F(n^{(\ell)})$, and
the fibres of $\Phi_V$ have dimension
$\sum_{\ell} n^{(\ell)}{}\tra D (m^{(\ell)}-n^{(\ell)})$. By induction, this yields
the bound
\begin{eqnarray*}
\dim(F(n)^g)&\leq& \sum_{\ell}\dim(F(n(\ell))+\dim(\ker(\Phi_V))\\
&=&\sum_{\ell}{}'\big(n^{(\ell)}{}\tra (D-I)n^{(\ell)}+1\big)
+ \sum_{\ell} n^{(\ell)}{}\tra D (m^{(\ell)}-n^{(\ell)}),
\end{eqnarray*}
where $\sum{}'$ signifies summation over all $\ell$ with $n^{(\ell)}\neq 0$.
Moreover, the centraliser $H\subset G(n)$ of $g$ is an open subset in
$$\left(\bigoplus_i\End(W_i)\right)^g\isom
\bigoplus_{\ell}\bigoplus_{i}\Hom(W_i^{(\ell)},K_i^{(\ell)})$$
and therefore has dimension $\dim(H)=\sum_{\ell}n^{(\ell)}{}\tra m^{(\ell)}$.
Connecting these pieces of information we obtain
\begin{eqnarray*}
\dim(G(n)F(n)^g)&\leq&\dim(F(n)^g)+\dim(G(n))-\dim(H)\\
&\leq &\sum_{\ell}{}'\big(n^{(\ell)}{}\tra (D-I)n^{(\ell)}+1\big)
+ \sum_{\ell} n^{(\ell)}{}\tra D (m^{(\ell)}-n^{(\ell)})\\
&&+n\tra n-\sum_\ell n^{(\ell)}{}\tra m^{(\ell)}.
\end{eqnarray*}
The difference of the last expression to the expected dimension of $F(n)$ is
$$
\Delta:=\Big[n\tra (D-I)n-n\tra n+1\Big]
-\sum_{\ell}{}'\Big[n^{(\ell)}{}\tra(D-I)m^{(\ell)}- n^{(\ell)}{}\tra n^{(\ell)}+1\Big].
$$
Note that the two bracketed expressions are not quite symmetric to
each other due to the presence of $m^{(\ell)}$ instead on
$n^{(\ell)}$. We can get rid of $n$ and $m^{(\ell)}$ due to the
relations
$$m^{(\ell)}=\sum_k n^{(k)}\min\{k,\ell\},\quad n=\sum_k n^{(k)}k,$$
and can rewrite the bound $\Delta$ in terms of the $n^{(k)}$ as follows:
$$
\Delta=
\sum_{\ell,k}n^{(k)}{}\tra (D-I)n^{(\ell)}\big(k\ell-\min\{k,\ell\}\big)
-\sum_{k,\ell}n^{(k)}{}\tra n^{(\ell)}k\ell+\sum_k{}'\Big(n^{(k)}{}\tra n^{(k)}-1\Big)+1
$$
Reorganise the sum in collecting those terms that contain $n^{(1)}$:
\begin{eqnarray*}%
\lefteqn{\Delta =1
+\Big[-1+2 n^{(1)}{}\tra \sum_{k\geq 2}\big((k-1)(D-2I)-I\big)n^{(k)}\Big]
+
\sum_{k\geq 2}{}'\big(n^{(k)}{}\tra n^{(k)}-1\big)}
\hspace{4em}
\\
&&+\sum_{k,\ell\geq 2}n^{(k)}{}\tra \Big((k\ell-\min\{k,\ell\})(D-2I)-\min\{k,\ell\} I\Big) n^{(\ell)}
\end{eqnarray*}
Here the second summand $[\ldots]$ appears only if $n^{(1)}\neq 0$. Note that there always
is at least one index $k\geq 2$ with $n^{(k)}\neq 2$, since $h\neq 0$. This shows
that all summands in the last expression for $\Delta$ are non-negative.

The minimal contribution of a non-zero vector $n^{(k)}$, $k\neq 2$, to $\Delta$ is
$$k((k-1)a -1)\left(\sum_i n^{(k)}_i\right)^2\geq 2.$$
Thus we always have $\Delta\geq 3$, and even better: $\Delta\geq 4$ in all cases except
$$a=2, n^{(2)}=(1), n^{(k)}=0\text{ for all }k\neq 2.$$
In this case $s=1$, $n=(2)$, and $d_{11}=4$, which is the same
exceptional case as before.

As in the semisimple case, the union of all sets $G(n)F(n)^g$, $g$
unipotent, is covered by a finite number of such sets. In fact, this
is even easier to see: up to conjugation there are only finitely
many different nilpotent elements $h$ and hence only finitely many
different subschemes $G(n)F(n)^{1+h}\subset Z$.

3.3. {\sl Case: $g\in G(n)\setminus\IC^*$ arbitrary.} Consider the multiplicative
Jordan decomposition $g=s u$, where $s$ is semisimple, $u$ is unipotent and $s$ and $u$
commute. Any endomorphism that commutes with $g$ also commutes with $s$ and $u$.
This implies that $F(n)^g\subset F(n)^s\cap F(n)^u$, so that the general case is
covered by 3.1. and 3.2. above.
\end{proof}

%%%%%%%%%%%%%%%%%%%%%%%%%%%%%%%%%%%%%%%%%%%%%%%%%%%%%%%%%%%%%%%%%%%%%%%%%%%%%%%%%%
\subsection{Return from the normal cone}\label{subsec:passagetonormalcone}
Let $v_0\in H^{\even}(X,\IZ)$ be a primitive Mukai vector satisfying $(*)$.
Let $v=mv_0$ for some multiplicity $m\in\IN$. We keep the notation introduced
earlier.

\begin{proposition}\label{prop:normalityofslice}---
Let $H$ be an arbitrary ample divisor. Let $E=\bigoplus_{i=1}^s E_i^{\oplus n_i}$
be a polystable sheaf whose stable direct summands $E_i$ satisfy the condition
\begin{equation}\label{eq:assumptiononvEi}
v(E_i)\in\IN v_0
\end{equation}
Consider a point $[q:\kh\to E]\in R^{ss}$ and a slice $S\subset
R^{ss}$ to the orbit of $[q]$ as above. Then $\ko_{S,[q]}$ is a
normal complete intersection domain of dimension
$$
\dim\Ext^1(E,E)-\dim\Ext^2(E,E)_0=
1+\sum_{i,j}n_i(\dim\Ext^1(E_i,E_j)-\delta_{ij})n_j,
$$
that has property $R_3$ in all cases except the following two:
\begin{enumerate}
\item $s=1$, $n_1=2$, $\dim\Ext^1(E_1,E_1)=4$,
\item $s=2$, $n_1=n_2=1$, $\dim\Ext^1(E_1,E_2)=2$.
\end{enumerate}
\end{proposition}

\begin{proof} Recall the notation introduced in sections
\ref{subsec:LocalDescription} and \ref{subsec:PassageToNormalCone}.
By Proposition \ref{prop:sympredfactorial}, $F=\mu^{-1}(0)=\Spec(\gr
A/J)$ is a normal complete intersection variety of dimension
\begin{eqnarray*}
\dim(F)&=&1+\sum_{i,j}n_i(\dim\Ext^1(E_i,E_j)-\delta_{ij})n_j\\
&=&\dim \Ext^1(E,E)-\dim \Ext^2(E,E)_0.
\end{eqnarray*}
Therefore, we must have equality at all places in inequality
\eqref{eq:dimestimate}.
Furthermore, since $F=\Spec(\gr A/J)$ is reduced and irreducible, the
equality of dimensions implies $J=\ini(I)$. It follows that
$$\gr (\widehat\ko_{S,[q]})=\gr(A/I)=\gr A/\ini(I)=\Gamma(F,\ko_F)$$
is a normal complete intersection. In particular, $\gr
(\widehat\ko_{S,[q]})$ is Cohen-Macaulay, hence satisfies $S_k$ for
all $k\in\IN$. Unless we are in the two exceptional cases, $\gr
(\widehat\ko_{S,[q]})$ is smooth in codimension 3. Now remark that
$\gr (\widehat\ko_{S,[q]})=\gr (\ko_{S,[q]})$
(\cite{AtiyahMacDonald}, 10.22) and then use the following
proposition which shows that $\ko_{S,[q]}$ itself is a normal
complete intersection which, unless we are in the two exceptional
cases, satisfies $R_3$.
\end{proof}

\begin{proposition}--- Let $B$ be a noetherian local ring
with maximal ideal $\gothm$ and residue field $B/\gothm\isom \IC$.
Let $\gr B$ denote the graded ring associated to the $\gothm$-adic
filtration of $B$. Then $\dim(B)=\dim(\gr B)$, and if $\gr B$ is an
integral domain or normal or a complete intersection then the same
is true for $B$. Moreover if $\gr B$ satisfies $R_k$ and $S_{k+1}$
for some $k\in\IN$ then $B$ satisfies $R_k$.
\end{proposition}

\begin{proof} The assertion about integrality and normality is Krull's theorem
(see \cite{Matsumura} (17.D) Thm 34). The assertions about complete
intersections and the property $R_k$ are due to Cavaliere and Niesi
(\cite{Cavaliere}, Theorems 3.4 and 3.13).
\end{proof}

\begin{lemma}\label{lem:Hgeneral}--- The assumption \eqref{eq:assumptiononvEi} in
Proposition \ref{prop:normalityofslice}
is satisfied in any of the following two situations:
\begin{enumerate}
\item $H$ is $v$--general.
\item $E=E_0^{\oplus m}$ for some stable sheaf $E_0$ with $v(E_0)=v_0$.
\end{enumerate}
The exceptions of Proposition \ref{prop:normalityofslice}
are met in case (B) only, i.\ e.\ if $\langle v_0,v_0\rangle =2$ and
$m=2$.
\end{lemma}

\begin{proof} Under the assumption that $H$ is $v$--general one has
$v(E_j)= r_jv_0$ for some $r_j\in\IN$ and all direct summands $E_j$
of $E$. Then $\dim \Ext^1(E_i,E_j)=r_ir_j\langle v_0,v_0\rangle \geq
2$. Thus Proposition \ref{prop:normalityofslice} applies.
%The exceptional cases are realised
%only if $\langle v_0,v_0\rangle=2$ and $r_1=1$ or $r_1=r_2$,
%respectively, i.e.\ in case (B).
\end{proof}

\begin{proposition}\label{prop:Risfactorial}---
1. Let $H$ be a $v$--general ample divisor. Then $R^{ss}$ is normal and
locally a complete intersection of dimension $\langle v,v\rangle+1+N^2$.
In case (C) it has property $R_3$ and hence is locally factorial.

2. Suppose that $E=E_0^{\oplus m}$ for some stable sheaf $E_0$ with
$v(E_0)=v_0$. Let $H$ be an arbitrary ample divisor. In case (C), there
is an open neighbourhood $U$ of $[E]\in M_v$ such that $\pi^{-1}(U)\subset
R^{ss}$ is locally factorial of dimension $\langle v,v\rangle+1+N^2$.
\end{proposition}

\begin{proof} 1. Let $[q]\in R^{ss}$ be a point with closed orbit, and let
$S\subset R^{ss}$ be a $\PAut(E)$--equivariant subscheme as in
subsection \ref{subsec:constructionofmoduli}. By Lemma
\ref{lem:Hgeneral} and Proposition \ref{prop:normalityofslice}, the
local ring $\ko_{S,[q]}$ is a normal complete intersection that has
property $R_3$ in case (C). But being normal or locally a complete
intersection or having property $R_k$ are open properties [EGA IV
19.3.3, 6.12.9]. Hence there is an open neighbourhood $U$ of $[q]$
in $S$ that is normal, locally a complete intersection, and has
property $R_3$ in case (C). The natural morphism $\LiePGl(N)\times
S\to R^{ss}$ is smooth. Therefore every closed orbit in $R^{ss}$ has
an open neighbourhood that is normal, locally a complete
intersection, and has property $R_3$ in case (C). Finally, every
$\LiePGl(N)$--orbit of $R^{ss}$ meets such an open neighbourhood. It
follows that $R^{ss}$ is normal, locally a complete intersection. In
case (C), $R^{ss}$ is regular in codimension 3 and hence locally
factorial due to the following theorem of Grothendieck.

2. The second assertion follows analogously.
\end{proof}

\begin{theorem} {\em (Grothendieck \cite{Grothendieck} Exp.\ XI Cor.\ 3.14)}
\label{thm:Grothendieck} --- Let $B$ be noetherian local ring. If
$B$ is a complete intersection and regular in codimension $\leq 3$,
then $B$ is factorial.
\end{theorem}

%%%%%%%%%%%%%%%%%%%%%%%%%%%%%%%%%%%%%%%%%%%%%%%%%%%%%%%%%%%%%%%%%%%%%%%%%%%%%%
%%%%%%%%%%%%%%%%%%%%%%%%%%%%%%%%%%%%%%%%%%%%%%%%%%%%%%%%%%%%%%%%%%%%%%%%%%%%%%

\section{A basic irreducibility result}\label{sec:irreducibility}

%%%%%%%%%%%%%%%%%%%%%%%%%%%%%%%%%%%%%%%%%%%%%%%%%%%%%%%%%%%%%%%%%%%%%%%%%%%%%%
%%%%%%%%%%%%%%%%%%%%%%%%%%%%%%%%%%%%%%%%%%%%%%%%%%%%%%%%%%%%%%%%%%%%%%%%%%%%%%

The following theorem generalises a beautiful result of Mukai \cite{Mukai}.

\begin{theorem}\label{thm:irreducibility}---
Let $X$ be a projective K3 or abelian surface with an ample
divisor $H$. Let $M_v$ be the moduli space of semistable sheaves associated to
a vector  $v\in H^{\even}(X,\IZ)$. Suppose that $Y\subset M_v$ is a connected
component parametrising stable sheaves only. Then $M_v=Y$.
\end{theorem}

\begin{proof} 1. Since all points in $Y$ correspond to stable sheaves, $Y$ is
smooth of expected dimension $\dim(Y)=2+\langle v,v\rangle$. Fix a point
$[F]\in Y$ and suppose that there is a point $[G]\in M_v\setminus Y$.
We shall exploit a beautiful old idea of Mukai \cite{Mukai}:
assume for a moment that there were a universal family
$\IF\in\text{Coh}(Y\times X)$. Let $p:Y\times X\to Y$ and $q:Y\times X\to X$
be the projections. We may then compare
the relative Ext-sheaves $\Ext_p^\bullet(q^*F,\IF)$ and $\Ext_p^\bullet(q^*G,\IF)$.
Since $F$ and $G$ are numerically equal on $X$, the same is true for the
classes of the Ext-sheaves according to the Grothendieck-Riemann-Roch theorem.
This will lead to a contradiction.

2. In general, there is no universal family, but the following construction will be
sufficient:

\begin{lemma}--- There is a smooth projective variety $Y'$ that parametrises
a family $\IF$ of stable sheaves on $X$ with Mukai vector $v$ such that the
classifying morphism $f:Y'\to Y$ is  surjective, generically finite, and
\'etale over a neighbourhood of $[F]$.
\end{lemma}

\begin{proof} Let $R':=Y\times_{M_v}R^{ss}$. Then $R'\to Y$ is a
$\LiePGl(N)$--principal fibre bundle, locally trivial in the \'etale topology.
Moreover, there is a universal epimorphism $\ko_{R'}\boxtimes \kh\to \IF'$.
We form the quotient $P:=(\IP^{N-1}\times R')\GIT \LieGl(N)$. Then $P$ is a
smooth projective variety, and the natural morphism $P\to Y$ is locally a product
in the \'etale topology with fibres isomorphic to $\IP^{N-1}$.
The center $\IC^*\subset\LieGl(N)$ acts trivially on the family
$\ko_{\IP^{N-1}}(-1)\boxtimes \IF'$. Therefore, this sheaf descends to a family
$\IF_P$ on $P\times X$. Let $L$ be a very ample line bundle on $P$. Choose a
linear subspace $Z\subset \IP(H^0(P,L))$ of codimension $N-1$ in such a way that
$Y':=Z\cap P$ is smooth and $f:Y'\to Y$ is \'etale over a neighbourhood of
$[F]$. Finally, let $\IF:=\IF_P|_{Y'\times X}$.
\end{proof}

3. Let $f:Y'\to Y$ and $\IF$ be chosen as in the lemma and let
$p:Y'\times X\to Y'$ and $q:Y'\times X\to X$ denote the two projections.
Moreover, let $f^{-1}([F])=\{p_1,\ldots,p_n\}$.

As $G$ represents a point in $M\setminus Y$ and hence is not isomorphic to
any of the stable sheaves $E$, $[E]\in Y$, one has $\Hom(G,E)=0=\Ext^2(G,E)$
for all $[E]\in Y$. It follows that $\Ext_p^0(q^*G,\IF)$
and $\Ext_p^2(q^*G,\IF)$ vanish and that $W:=\Ext^1_p(q^*G,\IF)$ is
a locally free sheaf on $Y'$ of rank $\langle v,v\rangle=\dim(Y)-2$.

If $G$ is replaced by $F$ the situation gets more complicated as the dimension
of the Ext-groups jumps on the fibre $T$. There is a complex of
locally free $\ko_{Y'}$--sheaves
\begin{equation}0\lra A^0\xra{\; \alpha\;}A^1\xra{\;\beta\;}A^2\lra 0\end{equation}
with the property that $\Ext_{p_{_S}}^i(t_{_X}^*q^*G,t_X^*\IF)\isom h^i(t^*(A^\bullet))$
for every base change
\begin{equation}\begin{array}{ccccc}
S\times X&\xra{\;t_X\;}&Y'\times X&\xra{\;q\;}&X\\
\scriptstyle{p_S}\Big\downarrow\phantom{\scriptstyle{p_S}}&&
\scriptstyle{p}\Big\downarrow\phantom{\scriptstyle{p}}\\
S&\xra{\;t\;}&Y'.
\end{array}
\end{equation}

\begin{lemma}--- The degeneracy locus of $\alpha$ and $\beta$ is the union of
the reduced points $p_1,\ldots,p_n$. Moreover, $\rk\alpha(p_i)=\rk A^0-1$ and
$\rk \beta(p_i)=\rk A^2-1$ for $i=1,\ldots,n$.
\end{lemma}

\begin{proof} For all $[E]\in Y$, $E\not\isom F$, one has $\Hom(F,E)=0=
\Ext^2(F,E)$. This implies that $\alpha$ and $\beta$ have maximal rank
on $Y'\setminus\{p_1,\ldots,p_n\}$. Moreover, $\Hom(F,F)=\IC=\Ext^2(F,F)$,
and this gives the second assertion of the lemma. It remains to show that
the degeneracy locus is reduced. Recall that tangent vectors in $T_{[F]}Y$
correspond bijectively to elements
$\gamma\in\Ext^1(F,F)$. Let $F_\gamma$ be the infinitesimal extension of $F$ over
$\Spec\;\IC[\varepsilon]$ corresponding to $\gamma$. The extension
\begin{equation}
0\lra F\xra{\;\varepsilon\;}F_\gamma\lra F\lra 0
\end{equation}
induces a long exact sequence
$$
\lra \Ext^i(F,F)\lra \Ext^i_{\Spec\,\IC[\varepsilon]}(F\tensor
\IC[\varepsilon],F_\gamma) \lra
\Ext^i(F,F)\xra{\;\partial\;}\Ext^{i+1}(F,F)\lra$$ where the
boundary operator is given by $\partial(e)=\gamma\cup e$. Now
$\gamma\cup -:\IC=\End(F,F)\to \Ext^1(F,F)$ is clearly injective,
and $\gamma\cup -:\Ext^1(F,F)\to \Ext^2(F,F)$ is surjective since
the symplectic form on $\Ext^1(F,F)$ is non-degenerate. It follows
that $\Ext^0(F,F)\isom \Ext^0_{\Spec\,\IC[\varepsilon]}(F\tensor
\IC[\varepsilon],F_\gamma)$ and
$\Ext^2_{\Spec\,\IC[\varepsilon]}(F\tensor
\IC[\varepsilon],F_\gamma)\isom \Ext^2(F,F)$. If the degeneracy
locus of $\alpha$ resp. $\beta$ were not reduced, the corresponding
$\Ext$ groups should be bigger than $\IC$ for at least one $\gamma$.
The calculation shows that this is not the case.
\end{proof}

4. Let $\sigma:Z\to Y$ denote the blow-up of $Y$ in $[F]$ with exceptional divisor $D$
and similarly $\varphi:Z'\to Y'$ the blow-up of $Y'$ in all points $p_i$ with
corresponding exceptional divisors $D_i$.
\begin{equation}
\begin{array}{ccccccc}
D&\subset&Z&\xla{\;g\;}&Z'&\supset&D_i\\
\Big\downarrow&&\scriptstyle{\sigma}\Big\downarrow\phantom{\scriptstyle{\sigma}}&
&\scriptstyle{\varphi}\Big\downarrow\phantom{\scriptstyle{\varphi}}&
&\scriptstyle{}\Big\downarrow\phantom{\scriptstyle{}}\\
{}[F]&\in&Y&\xla{\;f\;}&Y'&\ni&p_i
\end{array}
\end{equation}
According to the lemma, the degeneracy locus of both $\varphi^*(\alpha)$ and $\varphi^*(\beta)$ is
precisely the smooth divisor $D'=D_1\cup\ldots\cup D_n$. Therefore
these maps factor as follows:
\begin{equation} \varphi^*A^0\subset A'{}^0\xra{\;\alpha'\;}\varphi^*A^1\xra{\;\beta'\;}
A'{}^2\subset \varphi^*A^2,\end{equation}
where $A'{}^0$ and $A'{}^2$ are locally free, $\alpha'$ and $\beta'$ are
homomorphisms of maximal rank. Moreover, the line bundles $L:=\varphi^*A^2\big/A'{}^2$ and
$M:=A'{}^0\big/\varphi^*A^0$ on $D'$ are characterised by the canonical isomorphisms
$$L\tensor \ko_{D'}\isom \Ext_{D'}^2(q^*F,\ko_D\boxtimes\IF|_{f^{-1}([F])\times X})\isom \Ext^2(F,F)\tensor_\IC\ko_{D'}$$
and
$$\Tor_1^{\ko_{Z'}}(M,\ko_{D'})\isom
\Ext_{D'}^0(q^*F,\ko_D\boxtimes\IF|_{f^{-1}([F])\times X})
\isom \Hom(F,F)\tensor_\IC\ko_{D'},$$
implying
\begin{equation}\label{eq:linebundlesLandM}
L\isom \bigoplus_{i=1}^n\ko_{D_i}\quad\text{and}\quad
M\isom \bigoplus_{i=1}^n\ko_{D_i}(D_i).
\end{equation}

5. Let $W'$ denote the middle cohomology of the complex
$$
0\lra A'{}^0\xra{\;\alpha'\;} \varphi^*A^1\xra{\;\beta'\;} A'{}^2\lra 0.
$$
$W'$ is locally free of rank $\dim Y-2$. We obtain the following
equation of Chern classes in $H^*(Z',\IZ)$:
\begin{equation}\label{eq:chernidentity1}
\varphi^*c(A^1-A^0-A^2)=c(W'+M-L).%=c(W')\cdot c((\ko_D(D)-\ko_D)\boxtimes\ko_T(-1)).
\end{equation}
%where $x=1-\big(j_*((1-h)^{-1}\big)^2$ with $h=c_1(\ko_T)$ and $j:T'\to Z'$.
On the other hand, as $c(F)=c(G)$ in $H^*(X,\IZ)$, the
Grothendieck-Riemann-Roch Theorem yields the following identity in $H^*(Y',\IZ)$:
\begin{equation}\label{eq:chernidentity2}
c(A^1-A^0-A^2)=c(\Ext_p^\bullet(q^*F,\IF))=c(\Ext^\bullet_p(q^*G,\IF))=c(W).
\end{equation}
Combining \eqref{eq:chernidentity1} and \eqref{eq:chernidentity2}, we
conclude that
\begin{equation}\label{eq:chernidentity3}
c(W')=\varphi^*c(W)\cdot c(L-M)\in H^*(Z',\IZ)
\end{equation}
Moreover,
$$
c(L-M)=\prod_{i=1}^n\frac{c(\ko_{D_i})}{c(\ko_{D_i}(D_i))}=\prod_{i=1}^n
\frac1{c(\ko_{Z'}(-D_i))c(\ko_{Z'}(D_i))}=1+\sum_{k=1}^\infty
\sum_{i=1}^n D_i^{2k}.
$$
The product of any cohomology class in $H^*(Y',\IZ)$ of positive degree
with any of the classes $D_i$ is zero. It follows that
$$c_{2k}(W')=\varphi^*c_{2k}(W)+\sum_{i=1}^n D_i^{2k}\quad
\text{for all }k>0.$$ The key point now is that both $W$ and $W'$
are vector bundles of rank $\dim(Y)-2$, so that the Chern classes
$c_{\dim(Y)}(W)$ and $c_{\dim(Y)}(W')$ vanish (\emph{cf.}
\cite{Markman}, Lemma 4). We get the contradiction
$$0=\sum_{i=1}^n D_i^{\dim(Y)}=-n.$$
This finishes the proof of Theorem \ref{thm:irreducibility}.
\end{proof}

\begin{theorem}\label{thm:irred2}---
Let $v_0$ be a primitive Mukai vector satisfying condition
$(*)$ and $(**)$. Let $v=mv_0$ and let $H$ be a $v$--general ample divisor. Then
$M_v$ is a normal irreducible variety of dimension $2+\langle v,v\rangle$.
\end{theorem}

This theorem is due to Yoshioka \cite{Yoshtwisted} in the case of
torsion free sheaves. Using the local information obtained in
Proposition \ref{prop:Risfactorial}, the basic irreducibility result of
Theorem \ref{thm:irreducibility}, we can give a simple direct proof.

\begin{proof} By Proposition \ref{prop:Risfactorial}, $R^{ss}$ is normal.
As a GIT-quotient of a normal scheme, $M_v$ is also normal. If
$m=1$, all points in $M_v=M_{v_0}$ correspond to stable sheaves and
hence $M_{v}$ is smooth. By Theorem \ref{thm:irreducibility},
$M_{v_0}$ is irreducible. By $(**)$, $M_{v_0}$ is non-empty.

Assume now that $m\geq 2$ and that the assertion of the theorem has been proved
for all moduli spaces $M_{m'v_0}$, $1\leq m'<m$. For any decomposition
$m=m'+m''$ with $1\leq m'\leq m''$, consider the morphism
$$\varphi(m',m''):M_{m'v_0}\times M_{m''v_0}\lra M_{mv_0},\quad
([E'],[E''])\mapsto [E'\oplus E''],$$
and let $Y(m',m'')\subset M_v$ denote its image. The subschemes $Y(m',m'')$,
$1\leq m'\leq m''$,
are the irreducible components of the strictly semistable locus of $M_v$. Since
all $Y(m',m'')$ are irreducible by induction and intersect in the points of the
form $[E_0^{\oplus m}]$, $[E_0]\in M_{v_0}$, the strictly semistable locus is
connected. Since $M_v$ is normal, the connected components are irreducible.
In particular, there is exactly one component that meets the strictly semistable
locus. Theorem \ref{thm:irreducibility} excludes the possibility of a component
that does not meet the strictly semistable locus.
\end{proof}

%%%%%%%%%%%%%%%%%%%%%%%%%%%%%%%%%%%%%%%%%%%%%%%%%%%%%%%%%%%%%%%%%%%%%%%%%%%%%%%%%%
%%%%%%%%%%%%%%%%%%%%%%%%%%%%%%%%%%%%%%%%%%%%%%%%%%%%%%%%%%%%%%%%%%%%%%%%%%%%%%%%%%

\section{Factoriality of moduli spaces}\label{sec:Misfactorial}

%%%%%%%%%%%%%%%%%%%%%%%%%%%%%%%%%%%%%%%%%%%%%%%%%%%%%%%%%%%%%%%%%%%%%%%%%%%%%%%%%%
%%%%%%%%%%%%%%%%%%%%%%%%%%%%%%%%%%%%%%%%%%%%%%%%%%%%%%%%%%%%%%%%%%%%%%%%%%%%%%%%%%

\begin{proposition}\label{prop:Drezet}---
Let $v_0$ be a primitive Mukai vector satisfying $(*)$.
Let $v=mv_0$ for some $m\in\IN_0$. Assume that
\begin{itemize}
\item[--] either $E=E_0^{\oplus m}$, for some $E_0$ stable with $v(E_0)=v_0$, and
$H$ is arbitrary,
\item[--] or $E$ is arbitrary polystable with $v(E)=v$, and $H$ is $v$--general.
\end{itemize}
Assume further that case (C) applies.
Then $M_v$ is locally factorial at $[E]$ if and only if the isotropy subgroup
$\LiePGl(N)_{[q]}\isom \PAut(E)$ of any point $[q]$ in the closed orbit in
$\pi^{-1}([E])\subset R^{ss}$ acts trivially on the fibre $L([q])$ for
every $\LiePGl(N)$-linearised line bundle $L$ on an invariant open neighbourhood
of the orbit of $[q]$.
\end{proposition}

\begin{proof} This is Drezet's Th\'eor\`eme A \cite{Drezet}. In Drezet's
situation the Quot scheme $R^{ss}$ is smooth. However, all his arguments go
through under the weaker hypothesis that $R^{ss}$ is locally factorial in
a $\LiePGl(N)$--equivariant open neighbourhood of the closed orbit in the
fibre $\pi^{-1}([E])$. But this is true under the given hypothesis due to
Proposition \ref{prop:Risfactorial}
\end{proof}

\begin{corollary}\label{cor:mostdegenerate}---
Let $E_0$ be a stable sheaf with Mukai vector
$v(E_0)=v_0$ satisfying $(*)$ and assume that $v=mv_0$ satisfies (C).
Then $M_v$ is locally factorial at $[E_0^{\oplus m}]$.
\end{corollary}

\begin{proof} The isotropy subgroup of any point $[q]$ in the closed orbit
in $\pi^{-1}([E_0^{\oplus m}])\subset R^{ss}$ is isomorphic to
$\LiePGl(m)$ and therefore has no non-trivial characters. Hence the
action of $\LiePGl(m)$ on $L([q])$ is necessarily trivial (notations
as in Proposition \ref{prop:Drezet}).
\end{proof}

\begin{theorem}\label{thm:Misfactorial}---
Let $v_0$ be a primitive Mukai vector satisfying $(*)$ and $(**)$.
Assume that $v=mv_0$, $m\in\IN$, satisfies (C) and let $H$ be a $v$--general
ample divisor. Then $M_v$ is locally factorial.
\end{theorem}

\begin{proof}
Let $[E]\in M_v$ be an arbitrary point that is represented by the polystable
sheaf $E=\bigoplus_{i=1}^sE_i^{\bigoplus n_i}$, and let $[q:\kh\to E]$
be a point in the closed orbit in $\pi^{-1}([E])\subset R^{ss}$.
Since $H$ is $v$--general, the Mukai vectors of the stable direct summands
$E_i$ have the form
$$v(E_i)=m_i v_0,\quad m_i\in \IN, \quad \sum_{i=1}^s m_in_i=m.$$
We repeat the construction in section \ref{subsec:constructionofmoduli}
for each of the Mukai vectors $m_iv_0$, $i=1,\ldots,s$. Note that we can
choose a sufficiently large integer $k$ that works for all Mukai
vectors simultaneously. Let $P_i(z)= - m_i \langle v_0,v(\ko_X(- zH)\rangle$,
$N_i=P_i(k)$ and $\kh_i=\ko_{X}(-kH)^{\oplus N_i}$. Then $N=\sum_in_iN_i$ and
$\kh=\bigoplus_i\kh_i^{\oplus n_i}$. Moreover there are parameter spaces
$R_i^{ss}\subset\Quot_{X,H}(\kh_i,P_i)$ with $\LiePGl(N_i)$-actions and quotient
maps $\pi_i:R_i^{ss}\to M_{m_iv_0}$. Finally there is a canonical map
$$\Phi:\prod_i R_i^{ss}\lra R^{ss}, \quad
\big([\kh_i\to F_i]\big)_i\mapsto \big[\kh=\bigoplus_i\kh_i^{\oplus n_i}
\to \bigoplus_i F_i^{\oplus n_i}].$$
Let $Z$ denote the image of $\Phi$. It has the following properties:
\begin{itemize}
\item By Theorem \ref{thm:irred2}, the moduli spaces $M_{m_iv_0}$ are irreducible.
It follows that the schemes $R_i^{ss}$ and $Z$ are irreducible, too.
\item $Z$ contains the point $[q]$ and as well a point $[q':\kh\to E_0^{\oplus m}]$
for some stable sheaf $E_0$ with $v(E_0)=v_0$.
\item The group $G:=\big(\prod_i\LieGl(n_i)\big)/\IC^*\subset\LiePGl(N)$
fixes $Z$ pointwise. It equals the stabiliser subgroup of $[q]$ and is contained in
the stabiliser subgroup of $[q']$.
\end{itemize}
Now let $L$ be a $\LiePGl(N)$--linearised line bundle on $R^{ss}$.
The group $G$ acts on $L|_Z$ with a locally constant character, which must in fact
be constant, since $Z$ is connected. Moreover, the action is trivial at the point
$[q']$ according to the proof of Corollary \ref{cor:mostdegenerate}. Thus the
character is trivial everywhere on $Z$ and in particular at  $[q]$. According to
Drezet's criterion (Proposition \ref{prop:Drezet}), $M_v$ is locally factorial at
$[E]$.
\end{proof}

\begin{remark}
It is also known that the moduli space of semi-stable torsion free
sheaves on the projective plane is locally factorial by the work of
Drezet \cite{DrezetP2}. However it may be false for other surfaces
as has been observed by Le Potier: the moduli space
$M_{\IP^1\times\IP^1}(2,0,2)$ is not locally factorial at the point
represented by $\ko(1,-1)\oplus\ko(-1,1)$ (see \cite{Drezet}, p.
106).
\end{remark}

%%%%%%%%%%%%%%%%%%%%%%%%%%%%%%%%%%%%%%%%%%%%%%%%%%%%%%%%%%%%%%%%%%%%%%%%%%%%%%
%%%%%%%%%%%%%%%%%%%%%%%%%%%%%%%%%%%%%%%%%%%%%%%%%%%%%%%%%%%%%%%%%%%%%%%%%%%%%%

\section{Symplectic resolutions}\label{sec:SymplecticResolution}

%%%%%%%%%%%%%%%%%%%%%%%%%%%%%%%%%%%%%%%%%%%%%%%%%%%%%%%%%%%%%%%%%%%%%%%%%%%%%%
%%%%%%%%%%%%%%%%%%%%%%%%%%%%%%%%%%%%%%%%%%%%%%%%%%%%%%%%%%%%%%%%%%%%%%%%%%%%%%

Let $v_0$ be a primitive
Mukai vector satisfying $(*)$ and $(**)$. Let $v=mv_0$ and let $H$ be a
$v$--general divisor. Recall that the following three cases are possible:
\begin{itemize}
\item[(A)] $m=1$.
\item[(B)] $m=2$ and $\langle v_0,v_0\rangle=2$.
\item[(C)] $m\geq 3$, or $m=2$ and $\langle v_0,v_0\rangle \geq 4$.
\end{itemize}
In case (A) the moduli space $M_{v}$ consists only of stable sheaves.
It is irreducible and smooth of dimension $2+\langle v,v\rangle$.
Mukai \cite{Mukai} has defined a symplectic structure on $M_v$.

\begin{proposition}\label{prop:singularlocus}---
Assume that $m\geq 2$. The singular locus
$M_{v,\sing}$ of $M_v$ is non-empty and equals the semistable locus.
The irreducible components of $M_{v,\sing}$ correspond to integers
$m'$, $1\leq m'\leq m/2$, and have codimension
$2m'(m-m')\langle v_0,v_0\rangle-2$, respectively. In particular,
$\codim M_{v,\sing}=2$ in case (B) and $\geq 4$ in case (C).
\end{proposition}

\begin{proof} Recall the varieties
$Y(m',m'')$ introduced in the proof of Theorem \ref{thm:irred2}. The union of
the $Y(m',m'')$ is the strictly semistable locus. The maps
$$\varphi(m',m''):M_{m'v_0}\times M_{m''v_0}\to Y(m',m'')$$
are finite and surjective, hence
\begin{eqnarray*}%
\codim(Y(m',m''))&=&2+m^2\langle v_0,v_0\rangle-(2+m'{}^2\langle v_0,v_0\rangle)
-(2+m''{}^2\langle v_0,v_0\rangle)\\
&=&2m'm''\langle v_0,v_0\rangle-2.
\end{eqnarray*}
Clearly, the codimension 2 is attained only if $m'=m''=1$ and
$\langle v_0,v_0\rangle=2$, which is case (B). As $M_v$ is smooth in
all stable points, it remains to show that the strictly semistable
points are really singular. For this it suffices to show that $M_v$
is singular at a generic point $[E=E'\oplus E'']\in Y(m',m'')$,
where $E'$ and $E''$ are stable sheaves with  $v(E')=m'v_0$ and
$v(E'')=m''v_0$. In this case, $\PAut(E)\isom \IC^*$,
$\Ext^2(E,E)_0\isom\IC$, and the Kuranishi map $\Ext^2(E,E)_0\to
\IC[\Ext^1(E,E)]^{\wedge}$ is completely described by an invariant
function $f\in\IC[\Ext^1(E,E)]^{\wedge}$. It follows, that
$$\widehat\ko_{M_v,[E]}\isom
\big(\IC[\Ext^1(E,E)]^{\wedge}\big)^{\IC^*}\Big/(f).$$
Now $\IC^*$ acts on the four summands of
$$\Ext^1(E,E)=\Ext^1(E',E')\oplus \Ext^1(E',E'')
\oplus \Ext^1(E'',E')\oplus \Ext^1(E'',E'')$$
with weights $0$, $1$, $-1$, and $0$. It follows that
$$\Ext^1(E,E)\GIT\IC^*=\Ext^1(E',E')\times C\times \Ext^1(E'',E''),$$
where $C\subset M(d,\IC)$
is the cone of matrices of rank $\leq 1$ and
$$d=\dim \Ext^1(E',E'')=m'm''\langle v_0,v_0\rangle\geq 2.$$
Since the quotient of a singular local ring by a non-zero divisor
cannot become regular, $\widehat \ko_{M_v,[E]}$ is singular.
\end{proof}

\begin{theorem}--- Suppose that $v$ belongs to case (C). Then $M_v$ is a locally
factorial symplectic variety of dimension $2+\langle v,v\rangle$. The singular
locus is non-empty and has codimension 4. All singularities are symplectic, but
there is no open neighbourhood of a singular point in $M_v$ that admits a projective
symplectic resolution.
\end{theorem}

\begin{proof} We have already seen that $M_v$ is a locally factorial variety.
Mukai \cite{Mukai} constructed a non-degenerate 2-form on $M_v^s$. This form
is closed even if $M_v^s$ is not projective (\cite{HL} Prop.\ 10.3.2). By
Flenner's theorem \cite{Flenner} this form extends to any resolution of the
singularities of $M_v$. Hence the singularities are symplectic in the sense
of Beauville \cite{BeauvilleSymp}. Now let $[E]\in M_v$ be a singular point
and let $U\subset M_v$ be an open neighbourhood of $[E]$. A projective
symplectic resolution of $U$ is a projective resolution $\sigma:U'\to U$ of
the singularities of $U$ such that the restriction of the symplectic form on
$M_v^s$ to $U^{\reg}$ extends to a symplectic form on $U'$. In such a case
the morphism $\sigma$ would have to be semismall according to a result
of Kaledin, \cite{Kaledin} Lemma 2.11. As the singular locus of $U$ has
codimension $\geq 4$ according to Proposition \ref{prop:singularlocus}, the
exceptional locus of $\sigma$ has codimension $\geq 2$ in $U'$. On the other
hand $\ko_{M_v,[E]}$ is factorial by Theorem \ref{thm:Misfactorial}. This
implies that the exceptional locus must be a divisor (see \cite{Debarre}
no.\ 1.40 p.\ 28).
\end{proof}

\begin{remark}--- 1) The completion of a factorial local ring is not factorial in general.
The local rings of the moduli spaces of type (C) provide nice
examples of this phenomenon. Pushing the arguments in the previous
proof a bit further, one sees that
$$\widehat\ko_{M_v,[E]}\isom \IC[\Ext^1(E',E')\oplus \Ext^1(E'',E'')]^{\scriptstyle\wedge}\;\widehat\tensor\; B,$$
where $B$ is the completed coordinate ring of the cone $C_0\subset
C\subset M(d,\IC)$ of traceless matrices of rank $\leq 1$, with
$d\geq 4$. But $\widehat\ko_{M_v,[E]}$ cannot be factorial: the
vertex of $C_0$ is an isolated singularity of codimension $\geq 6$,
and there are two small symplectic resolutions
$T^*\IP(\Ext^1(E',E''))\rightarrow C_0\leftarrow
T^*\IP(\Ext^1(E'',E'))$. We see that in this case $\ko_{M,[E]}$ is
factorial due to Theorem \ref{thm:Misfactorial}, but
$\widehat\ko_{M,[E]}$ is not. Geometrically, what happens is this:
an irreducible Weil divisor becomes reducible after completion;
while the whole thing still is a Cartier divisor, some of its newly
acquired irreducible components need not be.

2) On the other hand, for polystable sheaves $E_0^{\oplus m}$ with
$E_0$ stable and $v(E_0)$ satisfying $(*)$, the completed local ring
$\widehat\ko_{M_{v},[E_0^{\oplus m}]}$ is factorial. In fact, the
proof of proposition \ref{prop:normalityofslice} shows that
$\widehat\ko_{S,[q]}$ is factorial. Moreover, the stabiliser is
isomorphic to $\LiePGl(m)$ hence has no non-trivial characters.
Under these conditions one can show that the invariant ring
$(\widehat\ko_{S,[q]})^{\LiePGl(m)}\simeq
\widehat\ko_{M_{v},[E_0^{\oplus m}]}$ is also factorial.
\end{remark}

%
% \subsection{Non-general ample divisors.} An dieser Stelle fehlt noch die
% Diskussion der Modulraume fr nicht-allgemeine ample Divisoren und die
% Aufloesung der Singularitaeten durch Uebergang zu benachbarten, allgemeinen
% Divisoren.
%
% %%%%%%%%%%%%%%%%%%%%%%%%%%%%%%%%%%%%%%%%%%%%%%%%%%%%%%%%%%%%%%%%%%%%%%%%%%%%%%%
% %%%%%%%%%%%%%%%%%%%%%%%%%%%%%%%%%%%%%%%%%%%%%%%%%%%%%%%%%%%%%%%%%%%%%%%%%%%%%%%
%
% \section{Case (B): 2-Factoriality}
%
% %%%%%%%%%%%%%%%%%%%%%%%%%%%%%%%%%%%%%%%%%%%%%%%%%%%%%%%%%%%%%%%%%%%%%%%%%%%%%%%
% %%%%%%%%%%%%%%%%%%%%%%%%%%%%%%%%%%%%%%%%%%%%%%%%%%%%%%%%%%%%%%%%%%%%%%%%%%%%%%%
%
% ******************************************************
%
% %%%%%%%%%%%%%%%%%%%%%%%%%%%%%%%%%%%%%%%%%%%%%%%%%%%%%%%%%%%%%%%%%%%%%%%%%%%%%%%
% %%%%%%%%%%%%%%%%%%%%%%%%%%%%%%%%%%%%%%%%%%%%%%%%%%%%%%%%%%%%%%%%%%%%%%%%%%%%%%%
%
% \section{Special cases}\label{sec:SpecialCases}
%
% Hier fehlt die Diskussion, sofern wir diese fuhren wollen, der Faelle
% $\langle v_0,v_0\rangle =-2$ und $=0$, sowie des Falls $v_0=(0,c_0,0)$.
%
%%%%%%%%%%%%%%%%%%%%%%%%%%%%%%%%%%%%%%%%%%%%%%%%%%%%%%%%%%%%%%%%%%%%%%%%%%%%%%%
%%%%%%%%%%%%%%%%%%%%%%%%%%%%%%%%%%%%%%%%%%%%%%%%%%%%%%%%%%%%%%%%%%%%%%%%%%%%%%%

\newpage
%%%%%%%%%%%%%%%%%%%%%%%%%%%%%%%%%%%%%%%%%%%%%%%%%%%%%%%%%%%%%%%%%%%%%%%%%%
%%%
%%%   Referenzliste
%%%
%%%%%%%%%%%%%%%%%%%%%%%%%%%%%%%%%%%%%%%%%%%%%%%%%%%%%%%%%%%%%%%%%%%%%%%%%%
\bibliographystyle{plain}

\begin{thebibliography}{10}
\bibitem{AtiyahMacDonald} M.~Atiyah, I.~Macdonald, Introduction to
Commutative Algebra, Addison-Wesley.

\bibitem{AK} A.~Altman, S.~Kleiman, Introduction to Grothendieck Duality
Theory. Lecture Notes in Mathematics 146. Springer Verlag 1970.

\bibitem{B} A.~Beauville, Vari\'et\'es K\"ahleriennes dont la premi\`ere
classe  de Chern est nulle.  J.\ Differential Geom. 18 (1983), 755 -- 782.

\bibitem{BeauvilleDonagi} A.~Beauville, Arnaud, R.~Donagi, La vari\'et\'e
des droites d'une hypersurface cubique de dimension 4.\ C.~R.~Acad.\ Sci.,
Paris, S\'er. I 301 (1985), 703 -- 706.

\bibitem{BeauvilleSymp} A.~Beauville, Symplectic singularities. Invent.\
math.\ 139 (2000), 541 -- 549.

\bibitem{Cavaliere} M.~P.~Cavaliere, G.~Niesi, On Serre's conditions in
the form ring of an ideal. J.~Math.\ Kyoto Univ. 21 (1981) 537 -- 546.

\bibitem{Debarre} O.~Debarre, Higher dimensional algebraic geometry.
Universitext Springer Verlag 2001.

\bibitem{WeilII} P.~Deligne, La conjecture de Weil, II.
  Inst. Hautes E'tudes Sci. Publ. Math. No. 52 (1980), 137--252.

\bibitem{DrezetP2} J.-M.~Drezet, Groupe de Picard des vari\'{e}t\'{e}s de modules de
faisceaux semi-stables sur $\IP^2_\IC$. Annales de l'institut
Fourier, 38 no. 3 (1988)

\bibitem{Drezet} J.-M.~Drezet, Points non factoriels des vari\'et\'es de
modules de faisceaux semi-stables sur une surface rationnelle. J.~reine
angew.\ Math. 413 (1991), 99 -- 126.

\bibitem{Flenner} H.~Flenner, Extendability of differential forms on
non-isolated singularities. Invent.\ math.\ 94 (1988), 317 -- 326.

\bibitem{Grothendieck} A.\ Grothendieck, Cohomologie locale des faisceaux
coh\'erents et Th\'eor\`emes de Lefschetz locaux et globaux. S\'eminaire de
G\'eom\'etrie Alg\'ebrique du Bois-Marie 1962 (SGA 2).
North-Holland Publishing Company Amsterdam 1968.

\bibitem{HL} D.~Huybrechts, M.~Lehn, The Geometry of Moduli Spaces of
Sheaves. Aspects of Mathematics E 31. Vieweg Verlag 1997.

\bibitem{Kaledin} D.~Kaledin, Symplectic resolutions from the Poisson point of view.
To appear in: J.~reine angewandte Math. math.AG/0310186.

\bibitem{Kiem} Y.~Kiem, On the existence of a symplectic desingularisation of some
moduli spaces of sheaves on a K3 surface. To appear in Comp.\ Math.

\bibitem{KiemChoy} Y.~Kiem, J.~Choy, Nonexistence of crepant resolution of moduli
space of sheaves on a K3 surface.

\bibitem{KL} D.~Kaledin, M.~Lehn, Local structure of hyperk\"ahler
singularities in O'Grady's examples. math.AG/0405575.

\bibitem{LS4} M.~Lehn, Ch.~Sorger, La singularit\'{e} de O'Grady.
math.AG/0504182.

\bibitem{Markman} E.~Markman, Generators of the cohomology ring of moduli spaces of sheaves
on symplectic surfaces.  J.\ reine angew.\ Math. 544 (2002), 61 --
82

\bibitem{Matsumura} H.~Matsumura, Commutative Algebra. W.\ A.\ Benjamin 1970.

\bibitem{Mukai} S.~Mukai, Symplectic structure of the moduli space of sheaves
on an abelian or K3 surface. Invent.\ math.\ 77 (1984), 101 -- 116.

\bibitem{Nakajima} H.~Nakajima, Instantons on ALE spaces, quiver
varieties, and Kac-Moody-algebras. Duke Math.~J.\ 76 (1994), 365 --
416.

\bibitem{OG1} K.~O'Grady, Desingularized moduli spaces of sheaves on a K3,
J.\ reine angew.\ Math. 512 (1999), 49 -- 117.

\bibitem{OG2} K.~O'Grady, A new six-dimensional irreducible symplectic variety.
J.\ Algebraic Geom. 12 (2003), 435 -- 505.

\bibitem{Rapagnetta} A.~Rapagnetta, Topological invariants of O'Grady's
six dimensional irreducible symplectic variety. math.AG/0406026.

\bibitem{Simpson} C.~Simpson, Moduli of representations of the fundamental group
of a smooth projective variety~I. Publ.\ Math.\ IHES 79 (1994), 47 -- 129.

\bibitem{Yoshioka} K.~Yoshioka, Moduli Spaces of stable sheaves on abelian
surfaces. Mathem.\ Ann.\ 321 (2001), 817 -- 884.

\bibitem{Yoshtwisted} K.~Yoshioka, Twisted stability and Fourier-Mukai Transform.
Comp.~Math. 138 (2003),\\ 261 -- 288.

\end{thebibliography}

\parindent0mm
\end{document}